\RequirePackage{fix-cm}

\documentclass[11pt,twoside,a4paper]{article}
\usepackage{cite}
\usepackage{datetime}\usdate
\usepackage[utf8]{inputenc}
\usepackage[english]{babel}
\usepackage{hyperref}

\usepackage{graphicx}
\usepackage{amsmath}
\usepackage{amsfonts}
\usepackage{algorithm}
\usepackage{algpseudocode} 
\usepackage{cite}
\usepackage{myartmods}

\algrenewcommand{\algorithmiccomment}[1]{\hfill[#1]}

\makeatletter

\def\blackslug{\hbox{\hskip 1pt \vrule width 4pt height 6pt depth 1.5pt
  \hskip 1pt}}
\let\qed\relax
\def\minimize#1{{\displaystyle\minim_{#1}}} 
\def\minim{\mathop{\operator@font minimize}}
\def\norm#1{\|#1\|}
\def\disp{\displaystyle}
\def\T{^T\!}
\def\inv{^{-1}}
\def\Tinv{^{-T}\!}
  \def\mtx#1#2{\renewcommand{\arraystretch}{1.2}%
      \left(\! \begin{array}{#1}#2\end{array}\! \right)}
\def\xhat{\skew{2.8}\widehat x}
\def\phat{\skew2\widehat p}
\def\ghat{\skew{4.3}\widehat g}
\def\text #1{\hbox{\quad#1\quad}}

\makeatother

\begin{document}

\title{ON EXACT LINESEARCH QUASI-NEWTON METHODS FOR MINIMIZING A
  QUADRATIC FUNCTION}

\author{Anders FORSGREN\thanks{\footTO} \and Tove
  ODLAND\addtocounter{footnote}{-1}\footnotemark}

\markboth{On exact linesearch quasi-Newton methods for minimizing a
  quadratic function}{}

\def\footTO{Optimization and Systems Theory, Department of
  Mathematics, KTH Royal Institute of Technology, SE-100 44 Stockholm,
  Sweden ({\tt andersf@kth.se,odland@kth.se}). Research partially supported by the Swedish Research Council (VR).}

\date{The final publication is available at Springer via
  \url{http://dx.doi.org/10.1007/s10589-017-9940-7}}

\maketitle

\begin{abstract}
  This paper concerns exact linesearch quasi-Newton methods for
  minimizing a quadratic function whose Hessian is positive
  definite. We show that by interpreting the method of conjugate
  gradients as a particular exact linesearch quasi-Newton method,
  necessary and sufficient conditions can be given for an exact
  linesearch quasi-Newton method to generate a search direction which
  is parallel to that of the method of conjugate gradients. 

  We also analyze update matrices and give a complete description of
  the rank-one update matrices that give search direction parallel to
  those of the method of conjugate gradients. In particular, we
  characterize the family of such symmetric rank-one update matrices
  that preserve positive definiteness of the quasi-Newton matrix. This
  is in contrast to the classical symmetric-rank-one update where
  there is no freedom in choosing the matrix, and positive
  definiteness cannot be preserved.

  The analysis is extended to search directions that are parallel to
  those of the preconditioned method of conjugate gradients in a
  straightforward manner.  


\end{abstract}

\section{Introduction}
In this paper we study the behavior of quasi-Newton methods (QN) on an
unconstrained quadratic problem of the form
\begin{equation}\label{qp}
\minimize{x \in \mathbb{R}^n} \quad \frac{1}{2}x^T Hx+c^T x,\tag{QP}
\end{equation}
where $H=H^T \succ 0$. Solving (\ref{qp}) is equivalent to solving a
symmetric system of linear equations $Hx+c=0$. In particular, our
concern is to give conditions under which a quasi-Newton method
utilizing exact linesearch generates search directions that are
parallel to those of the method of conjugate gradients (CG). As exact
linesearch is considered, parallel search directions imply identical
iterates. At iteration $k$, the $x$-iterate and the gradient $Hx+c$
are denoted by $x_k$ and $g_k$ respectively. In a quasi-Newton method,
the search direction $p_k$ is computed from $B_k p_k=-g_k$, where
$B_k$ is nonsingular.

We give necessary and sufficient conditions on a QN-method for this
equivalence with CG on (\ref{qp}). This is not the first time
necessary and sufficient conditions are given.  In \cite[Theorem
2.2]{KON98}, a necessary and sufficient condition is given, which is
based on projections from iterations $0$, $1$, \dots, $k-1$, allowing
also the preconditioned setting to be considered. In contrast, we
interpret the method of conjugate gradients as a particular
quasi-Newton method and base the necessary and sufficient conditions
on this observation. The result we give is thus directly based on the
projection given by the method of conjugate gradients, i.e., based on
quantities from iteration $k-1$ and $k$ involving one projection only.

If considering \emph{update matrices} $U_k$ defined by
$U_k=B_k-B_{k-1}$, it is well-known that, on (\ref{qp}), QN using
exact linesearch and an update scheme in the one-parameter Broyden
family generates identical iterates to those generated by CG, see,
e.g., \cite{fletcherpractical,huang,Nazareth}. The unique rank-1
update matrix in the Broyden family is usually referred to as the SR1
update matrix, and it is determined entirely by the so-called secant
condition. As a result of our equivalence result, we show that the
symmetric rank-1 update matrices that give parallel search directions
to CG are given by the family of update matrices on the form
\begin{equation}\label{eqn-rank1intro}
U_k = -
\frac1{
(\gamma_k-1) p_{k-1}^T g_{k-1}}
(\gamma_{k}g_k - g_{k-1})(\gamma_{k}g_k-g_{k-1})^T,
\end{equation}
where $\gamma_k$ is a free parameter. The free parameter can be seen
as a relaxation of the secant condition, as the SR1 update matrix is
the only matrix in our parameterized rank-1 family which satisfies
this condition. We show how to choose the parameter so that positive
definiteness of the quasi-Newton matrix is preserved.

To simplify the exposition, we discuss equivalence to CG, which
corresponds to the initial Hessian approximation being the identity
matrix in the quasi-Newton method in our analysis. We then give the
corresponding results in the preconditioned setting, which corresponds
to an arbitrary positive definite and symmetric initial Hessian
approximation. For the rank-1 case, the family of symmetric update
matrices take the form (\ref{eqn-rank1intro}) also in the
preconditioned setting.

In Section~\ref{background}, we make a brief introduction to CG and
QN. In Section~\ref{sec-results}, we present our results which include
necessary and sufficient conditions on QN such that CG and QN generate
parallel search directions. These results are specialized to update
matrices in Section~\ref{sec-update}. In particular, in
Section~\ref{sec-rank1}, we give the results on symmetric rank-1
update matrices. Section~\ref{sec-approx} contains a discussion on how
the results would apply if the inverse of the Hessian was updated
instead of the Hessian itself. In Section~\ref{sec-precond}, the
corresponding results in the preconditioned setting are
stated. Finally, in Section~\ref{Conclusion} we make some concluding
remarks.

\section{Background}\label{background}

For solving (\ref{qp}), we consider linesearch methods on the
following form. At iteration $k$, a search direction $p_k$ is
computed. The $x$-iterate and the gradient are updated as
\[
x_{k+1}=x_k+\theta_k p_k, \quad g_{k+1}=g_k+\theta_k Hp_k,
\text{for}
\theta_k=-\frac{g_k^T p_k}{p_k^T Hp_k}.
\]
The choice of steplength $\theta_k$ corresponds to exact linesearch,
i.e., given a search direction $p_k$ the steplength gives the exact
minimizer along $p_k$. This is a natural choice for (\ref{qp}), as it
can be done explicitly. For a given initial point $x_0$, the iteration
process is terminated at an iteration $r$ if $g_{r}=0$, in which case
$x_{r}$ is given as the optimal solution to (\ref{qp}) or equivalently
as the unique solution to $Hx+c=0$. The method is summarized in
Algorithm~\ref{alg-linesearch}.

\begin{algorithm}
\caption{\label{alg-linesearch}An exact linesearch method for solving $Hx+c=0$.}
\begin{algorithmic}[0]
\State $k\gets 0$;
\quad $x_k \gets$ initial point; \quad $g_k \gets H x_k + c$;
\While{$\norm{g_k} \ne 0$}
\State  $p_k \gets$ search direction;
\State $\disp\theta_k \gets -\frac{g_k\T p_k}{p_k\T
  H p_k}$;
\State  $x_{k+1} \gets x_k+\theta_k p_k$;
\quad  $g_{k+1} \gets g_k+\theta_k H p_k$;
\State  $k \gets k + 1$;
\EndWhile
\end{algorithmic}
\end{algorithm}

The particular linesearch method is defined by the way the search
direction $p_k$ is obtained in each iteration $k$. Our model method is
the method of conjugate gradients, CG, by Hestenes and Stiefel
\cite{HestenesStiefel}. There are different varieties of CG, which are
equivalent on (\ref{qp}). The variety we describe is referred to as
the Fletcher-Reeves method of conjugate gradients, as stated in the
following definition.

\begin{definition}[The method of conjugate gradients (CG)]\label{def-cg}
The \emph{method of conjugate gradients}, CG, is the linesearch
  method of the form given by Algorithm~$\ref{alg-linesearch}$ in
  which the search direction $p_k$ is given by $p_k^{CG}$, with
\begin{equation}\tag{CG}
p_k^{CG}=\begin{cases}
-g_0 & \mbox{if } k=0, \\
-g_k+\frac{g_k^Tg_k}{g_{k-1}^Tg_{k-1}}p_{k-1}^{CG} & \mbox{if } k\ge 1.
\end{cases}
\end{equation}
\end{definition}
For CG it holds that, for all $k$, $g_k^Tg_i=0$, $i=0, \dots, k-1$, so
the method terminates with $g_r=0$ for some $r$, $r\le n$, and $x_r$
solves (\ref{qp}).  In addition, it holds that
$\{p_k^{CG}\}_{k=0}^{r-1}$ are mutually conjugate with respect to
$H$. For an introduction to CG, see, e.g.,
\cite{cgwopain,saad,demmel}. In \cite{fletcherreeves}, CG is extended
to general unconstrained problems. The reason for CG being our model
method is that it requires one matrix-vector product $Hp_k$ per
iteration, and it terminates in $r$ iterations, with $r\le n$.

Next we define what we will refer to as a quasi-Newton method, QN. 

\begin{definition}[Quasi-Newton method (QN)]
A \emph{quasi-Newton method}, QN, is a linesearch method of the form
  given by Algorithm~\ref{alg-linesearch} in which the search
  direction $p_k$ is given by
\begin{equation}\tag{QN} 
B_k p_k=-g_k,
\end{equation}
where the matrix $B_k$ is assumed nonsingular.
\end{definition}

Quasi-Newton methods were first suggested by Davidon, see
\cite{davidon}, and later modified and formalized by Fletcher and
Powell, see \cite{fletcherpowell}. For an introduction to QN-methods,
see, e.g., \cite[Chapter 4]{practicalopt}.

Our interest is now to set up conditions on $B_k$ such that $p_k$ and
$p_k^{CG}$ are parallel for all $k$, so that QN also terminates in $r$
iterations. In \cite{ForsgrenOdland}, we derived such conditions based
on a sufficient condition to obtain mutually conjugate search
directions. Here, we give a direct necessary and sufficient condition
based on $p_k^{CG}$ only.

The results of the paper are derived with (CG) as the model method,
which corresponds to $B_0=I$ in (QN) giving $p_0=p_0^{CG}$. It is also
of interest to consider the case when a symmetric positive definite
matrix $M$ is given for which a \emph{preconditioned} method of
conjugate gradients is defined. This corresponds to $B_0=M$ in (QN)
giving the initial search directions identical. To simplify the
exposition, we derive the results for the unpreconditioned case given
by (CG) and give the corresponding results for the preconditioned
setting in Section~\ref{sec-precond}.

\section{Necessary and sufficient conditions for QN}\label{sec-results}

In this section we give precise conditions on $B_k$ such that $p_k$ is
parallel to $p_k^{CG}$. The main benefit of the conditions compared to
previous work is that our result is based on the single iteration
$k$. The dependence on the previous iterates is contained in the
search direction $p_{k-1}$, and there is no need to check any
condition for all the previous iterates.

In the following proposition, we give a necessary and sufficient
condition on $B_k$ at a particular iteration $k$ to give a search
direction $p_k$ such that $p_k=\delta_k p_k^{CG}$ for a scalar
$\delta_k$. We assume that each previous search direction $p_i$ has
been parallel to the corresponding search direction of CG, $p_i^{CG}$,
so that QN and CG have generated the same iterate $x_k$.

\begin{proposition}\label{prop-iffB}
  Consider iteration $k$ of the exact linesearch method of
  Algorithm~$\ref{alg-linesearch}$, where $1\le k < r$. Assume that
  $p_i=\delta_i p_i^{CG}$ with $\delta_i\ne 0$ for $i=0,\dots,k-1$,
  where $p_i^{CG}$, $i=0,\dots,k-1$, are the search directions of the
  method of conjugate gradients, as stated in
  Definition~$\ref{def-cg}$. Let $A_k$ be
  defined as
\begin{equation}\label{eqn-Ak}
A_k = I-\frac{1}{g_{k-1}^T p_{k-1}} p_{k-1} g_k^T.
\end{equation}
Then,
\begin{equation}\label{eqn-Akinv}
A_k\inv = I+\frac{1}{g_{k-1}^T p_{k-1}} p_{k-1} g_k^T,
\end{equation}
and it holds that $A_k p_k^{CG}=- g_k$. In addition, if $p_k$ is
given by $B_k p_k=-g_k$ with $B_k$ nonsingular, then, for any nonzero
scalar $\delta_k$, it holds that $p_k=\delta_k p_k^{CG}$ if and only
if
\begin{equation}\label{eqn-iffB}
  B_k A_k\inv g_k = 
\frac1{\delta_k} g_k,
\end{equation}
or equivalently if and only if
\begin{equation}\label{eqn-iffW}
B_k = A_k^T W_k A_k, \text{with} W_k g_k = \frac1{\delta_k} g_k, \text{for $W_k$ nonsingular.}
\end{equation}
Finally, it holds that $B_k\succ 0$ if and only if $W_k\succ 0$.
\end{proposition}

\begin{proof}
We have
\begin{equation}\label{eqn-pkCG1}
p_k^{CG} =- g_k +\frac{g_k^T g_k}{g_{k-1}^Tg_{k-1}}p_{k-1}^{CG}
= -\big(I+\frac{1}{g_{k-1}^T p_{k-1}^{CG}}p_{k-1}^{CG} g_k^T\big)g_k.
\end{equation}
Therefore, since $p_{k-1}=\delta_{k-1}p_{k-1}^{CG}$, with
$\delta_{k-1}\ne 0$, (\ref{eqn-pkCG1}) gives
\begin{equation}\label{eqn-pkCG2}
p_k^{CG} = -\big(I+\frac{1}{g_{k-1}^T p_{k-1}}p_{k-1} g_k^T\big)g_k=
-A_k\inv g_k,
\end{equation}
with $A_k\inv$ given by (\ref{eqn-Akinv}). Since $g_k^T p_{k-1}=0$,
multiplication of $A_k$ of (\ref{eqn-Ak}) by $A_k\inv$ of
(\ref{eqn-Akinv}) gives $A_k A_k\inv=I$, so that the stated $A_k$ is
nonsingular with corresponding inverse $A_k\inv$. Therefore,
(\ref{eqn-pkCG2}) gives $A_k p_k^{CG}=-g_k$.

Since $B_k$ is assumed nonsingular and $\delta_k\ne 0$, it holds that
$p_k=\delta_k p_k^{CG}$ if and only if $B_k (-\delta_k A_k\inv g_k) =
-g_k$, which is equivalent to (\ref{eqn-iffB}). Since $g_k^T
p_{k-1}=0$, we obtain $A_k\Tinv g_k=g_k$, so that (\ref{eqn-iffB}) is
equivalent to 
\[
A_k\Tinv B_k A_k\inv g_k=\frac1{\delta_k} g_k,
\]
which in turn is equivalent to (\ref{eqn-iffW}). The final result on
positive definiteness follows from the nonsingularity of $A_k$ by
Sylvester's law of inertia, see, e.g., \cite[Theorem
8.1.17]{golubvanloan}.  \qed \end{proof}

The necessary and sufficient conditions of Proposition~\ref{prop-iffB}
give a straightforward way to check if a matrix $B_k$ is such that the
corresponding QN-method and CG will generate parallel search
directions. The scaling of $p_k^{CG}$ has a special role in our
analysis and we relate $p_k$ to $p_k^{CG}$ by a scalar $\delta_k$. The
observation that $p_k^{CG}$ may be written as $p_k^{CG}=-A_k\inv g_k$
for a nonsingular $A_k$ has been made in \cite[Example 2]{KON98}, but
the equivalence result of \cite{KON98} concerns $p_k$ without relating
to the scaling of $p_k^{CG}$ explicitly. Therefore, the condition of
\cite{KON98} involves projections on all previous iterations
$0,1,\dots,k-1$, not one single projection as we obtain. In addition,
since there is no relationship to a particular scaling, there is no
parameter corresponding to our $\delta_k$. Since such a parameter is
vital for deriving later results in our paper, in particular when
characterizing symmetric rank-one updates, we cannot apply the
equivalence result of \cite{KON98} directly. A difference in
\cite{KON98} is that they consider matrices $N_k$ that approximate
$H\inv$ rather than matrices $B_k$ that approximate $H$. This is not a
major difference, we discuss these issues in Section~\ref{sec-approx}.

Note that it is not necessary to make $B_k-I$ increase in rank. In
particular, $B_k = A_k^T A_k$, corresponding to $W_k=I$ in
Proposition~\ref{prop-iffB}, is a positive-definite symmetric matrix
for which $B_k p_k=-g_k$ gives $p_k=p_k^{CG}$.

We also note that the characterization of $B_k$ does only depend on
information from iteration $k$ and $k-1$, since it directly inherits
the properties of the method of conjugate gradients. In addition, the
characterization of Proposition~\ref{prop-iffB} only depends directly
on quantities computed by the quasi-Newton method, the scaling of the
method of conjugate gradients is not needed.

\section{Results on update matrices}\label{sec-update}

In the previous section we gave results on $B_k$ for a particular
iteration $k$ without directly relating to any other $B_i$, $i\ne
k$. It is often the case that $B_k$ is defined in terms of the
previous matrix $B_{k-1}$ and an update matrix $U_k$ such that
$B_k=B_{k-1}+U_k$, and that conditions are put on $U_k$. We have in
mind a setting where information from the generated gradients is used,
so that $B_{k-1}$ may be expressed as $B_{k-1}=I+V_k$, with
$\mathcal{R}(V_k) \subseteq span\{g_0, \dots,g_{k-1}\}$. As we have in
mind such a setting where in addition $B_k$ is symmetric, we make the
assumption $B_{k-1}g_k=g_k$.

Proposition~\ref{prop-iffB} can then be applied in a straightforward
manner to give conditions on $U_k$ such that $p_k=\delta_k
p_k^{CG}$. Note that there is a one-to-one correspondence between
$U_k$ and $B_k$ given $B_{k-1}$.

\begin{proposition}\label{prop-iffU}
  Consider iteration $k$ of the exact linesearch method of
  Algorithm~$\ref{alg-linesearch}$, where $1\le k < r$. Assume that
  $p_i=\delta_i p_i^{CG}$ with $\delta_i\ne 0$ for $i=0,\dots,k-1$,
  where $p_i^{CG}$, $i=0,\dots,k-1$, are the search directions of the
  method of conjugate gradients, as stated in
  Definition~$\ref{def-cg}$.  Let $B_{k-1}$ be a nonsingular matrix
  such that $B_{k-1} p_{k-1}=-g_{k-1}$ and $B_{k-1}g_k = g_k$. Let
  $U_{k}=B_k - B_{k-1}$ and assume that $B_k$ and $p_k$ satisfy
  $B_kp_k=-g_k$, with $B_k$ nonsingular. Then, for any nonzero scalar
  $\delta_k$, it holds that $p_k=\delta_k p_k^{CG}$ if and only if
\begin{equation}\label{condU}
U_{k} \big( g_k+\frac{g_k^Tg_k}{p_{k-1}^T
  g_{k-1}}p_{k-1} \big) = (\frac1{\delta_k} -1 ) g_k + \frac{ g_k^T
  g_k}{p_{k-1}^T g_{k-1}} g_{k-1}.
\end{equation}
\end{proposition}

\begin{proof}
  By assumption, $B_k$ is nonsingular so for $B_k=B_{k-1}+U_k$,
  Proposition~\ref{prop-iffB} gives $p_k=\delta_k p_k^{CG}$ if and
  only if
\begin{align*}
U_k \big( g_k+\frac{g_k^Tg_k}{p_{k-1}^T
  g_{k-1}}p_{k-1} \big) &
=\frac1{\delta_k}g_k-B_{k-1}\big( g_k+\frac{g_k^Tg_k}{p_{k-1}^T
  g_{k-1}}p_{k-1} \big) \\
&=\frac1{\delta_k}g_k- g_k+\frac{g_k^Tg_k}{p_{k-1}^T g_{k-1}}g_{k-1},
\end{align*}
since $p_{k-1}$ is computed from $B_{k-1}p_{k-1}=-g_{k-1}$ and it is
assumed that $B_{k-1} g_k = g_k$, so the statement of the proposition
follows.
\qed \end{proof}

Note that in the right-hand side of (\ref{condU}) in
Proposition~\ref{prop-iffU}, the component along $g_{k-1}$ is nonzero
and independent of $\delta_k$. The component along $g_k$ is zero for
$\delta_k=1$, i.e., when $p_k=p_k^{CG}$.

\subsection{Results on symmetric rank-one update matrices}\label{sec-rank1}

Next we consider the case when $U_k$ is a symmetric matrix of rank
one. It is well known that the secant condition gives a unique update
referred to as SR1, see, e.g., \cite[Chapter 9]{luenberger}. The
secant condition and SR1 will be discussed later in this
section. Using Proposition~\ref{prop-iffU} we can give a different
result concerning the case when the update matrix $U_k$ is a symmetric
matrix of rank one. In particular, we show that the family of rank-1
update matrices can be parameterized by a free parameter and that the
matrix is unique for a fixed value of the parameter. This
parametrization allows positive definiteness of the quasi-Newton
matrix to be preserved.

The situation can be considered in two ways. First, for any given
value of the scalar $\gamma_k$, except for three distinct values,
there is a symmetric rank-1 update matrix $U_k$ of the form
\begin{equation}\label{eqn-Uk}
U_k = -
\frac1{
(\gamma_k-1) p_{k-1}^T g_{k-1}}
(\gamma_{k}g_k - g_{k-1})(\gamma_{k}g_k-g_{k-1})^T,
\end{equation}
for which $p_k=\delta_k(\gamma_k) p_k^{CG}$, where $\delta_k(\cdot)$
is a real-valued function. Second, if $p_k=\delta_k p_k^{CG}$ is
required for any given value of the scalar $\delta_k$, except for
three distinct values, and $U_k$ is symmetric and of rank one, $U_k$
must take the form (\ref{eqn-Uk}), with $\gamma_k=\gamma_k(\delta_k)$,
where $\gamma_k(\cdot)$ is the inverse function of
$\delta_k(\cdot)$. Consequently, except for three distinct values,
there is a one-to-one correspondence between $\delta_k$ such that
$p_k=\delta_k p_k^{CG}$ and $\gamma_k$ of the symmetric rank-1 update
matrix $U_k$ of (\ref{eqn-Uk}).

The functions $\delta_k(\cdot)$ and $\gamma_k(\cdot)$ are defined in
the following lemma.

\begin{lemma}\label{lem-gammadelta}
  Consider iteration $k$ of the exact linesearch method of
  Algorithm~$\ref{alg-linesearch}$, where $1\le k < r$. Assume that
  $p_i=\delta_i p_i^{CG}$ with $\delta_i\ne 0$ for $i=0,\dots,k-1$,
  where $p_i^{CG}$, $i=0,\dots,k-1$, are the search directions of the
  method of conjugate gradients, as stated in
  Definition~$\ref{def-cg}$. Let $\hat\gamma_k={p_{k-1}^T
    g_{k-1}}/{g_k^T g_k}$. For $\delta_k\ne 0$ and
  $\gamma_k\ne\hat\gamma_k$, let the functions $\gamma_k(\delta_k)$
  and $\delta_k(\gamma_k)$ be defined by
\[
\gamma_k(\delta_k) =- \frac{p_{k-1}^T g_{k-1}}{g_k^T g_k}
\left(\frac1{\delta_k}-1\right), \quad
\delta_k(\gamma_k) =\frac1{1 - \gamma_k\frac{g_k^T
    g_k}{p_{k-1}^T g_{k-1}}}.
\]
Then, the functions $\gamma_k(\cdot)$ and $\delta_k(\cdot)$ are
inverses to each other.
\end{lemma}

We now characterize the symmetric rank-one update matrices that give
search directions which are parallel to those of the method of
conjugate gradients. In addition, we give conditions for preserving
positive definiteness and a hereditary result.

\begin{proposition}\label{prop-rank1}
  Consider iteration $k$ of the exact linesearch method of
  Algorithm~$\ref{alg-linesearch}$, where $1\le k < r$. Assume that
  $p_i=\delta_i p_i^{CG}$ with $\delta_i\ne 0$ for $i=0,\dots,k-1$,
  where $p_i^{CG}$, $i=0,\dots,k-1$, are the search directions of the
  method of conjugate gradients, as stated in
  Definition~$\ref{def-cg}$. Let $B_k$ and $p_k$ satisfy
  $B_kp_k=-g_k$, and let $B_{k-1}$ be a nonsingular matrix such that
  $B_{k-1} p_{k-1}=-g_{k-1}$ and $B_{k-1}g_k = g_k$. In addition, let
  $\gamma_k(\cdot)$, $\delta_k(\cdot)$ and $\hat\gamma_k$ be given by
  Lemma~$\ref{lem-gammadelta}$.

  For any scalar $\gamma_k$, except $\gamma_k=0$,
  $\gamma_k=\hat\gamma_k$ and $\gamma_k=1$, let $B_k$ be defined by
\begin{equation}\label{Urank1}
B_{k} = B_{k-1} -
\frac1{
(\gamma_k-1) p_{k-1}^T g_{k-1}}
(\gamma_{k}g_k - g_{k-1})(\gamma_{k}g_k-g_{k-1})^T.
\end{equation}
Then, $B_k$ is nonsingular and $p_k=\delta_k p_k^{CG}$ for
$\delta_k=\delta_k(\gamma_k)$.

Conversely, for any scalar $\delta_k$, except $\delta_k=0$,
$\delta_k=\delta_k(1)$ and $\delta_k=1$, assume that $p_k=\delta_k
p_k^{CG}$ and assume that $B_k-B_{k-1}$ is symmetric and of rank
one. Then, $B_k$ is a nonsingular matrix given by $(\ref{Urank1})$ for
$\gamma_k=\gamma_k(\delta_k)$.

If, in addition, $B_{k-1}=B_{k-1}^T\succ0$, then $B_k$ defined by
$(\ref{Urank1})$ satisfies $B_k\succ 0$ if and only if $\gamma_k>1$ or
$\hat\gamma_k<\gamma_k<0$, or equivalently if and only if
$\gamma_k=\gamma_k(\delta_k)$ for $0<\delta_k<\delta_k(1)$ or $\delta_k>1$.

Finally, if $B_i p_i=-g_i$, $i=0,\dots,k$, with $B_0=I$ and if, for
$i=1,\dots,k$, $B_{i-1}$ is updated to $B_i$ according to
$(\ref{Urank1})$ for $\gamma_i$ such that $\gamma_i\ne0$,
$\gamma_i\ne\hat\gamma_i$ and $\gamma_i\ne1$, then
\begin{equation}\label{eqn-Bkpi}
B_k p_i = \frac{\gamma_{i+1}\theta_i}{\gamma_{i+1}-1}H p_i, \quad
i=0,\dots,k-1.
\end{equation}
\end{proposition}

\begin{proof}
  Let $U_k=B_k-B_{k-1}$. If $U_k$ is symmetric and of rank one, we may
  write $U_k=\beta_k u_k u_k^T$, where $\beta_k$ is a scalar and $u_k$
  is a vector in $\mathbb{R}^n$, both to be determined. If $B_k$ is
  nonsingular and $\delta_k\ne 0$, Proposition~\ref{prop-iffU} shows
  that $p_k=\delta_k p_k^{CG}$ if and only if
\begin{equation}\label{eqn-beta1}
\beta_k u_k u_k^TA_k^{-1} g_k = (\frac1{\delta_k} -1 ) g_k + \frac{ g_k^T
  g_k}{p_{k-1}^T g_{k-1}} g_{k-1}=
- \frac{g_k^T g_k}{p_{k-1}^T g_{k-1}}(\gamma_{k} g_k-g_{k-1}),
\end{equation}
with $\gamma_k=\gamma_k(\delta_k)$ given by
Lemma~\ref{lem-gammadelta}. Throughout the proof, assume that
$\delta_k\not\in\{0,\delta_k(1),1\}$ and
$\gamma_k\not\in\{0,\hat\gamma_k,1\}$, which is assumed in the
statement of the Proposition. Then, Lemma~\ref{lem-gammadelta} shows
that there is a one-to-one correspondence between $\delta_k$ and
$\gamma_k$. Hence, (\ref{eqn-beta1}) may be considered for either
$\delta_k$ or $\gamma_k$. We choose $\gamma_k$ for ease of notation.

We first assume that $B_k$ is nonsingular, and verify that this is the
case later in the proof. For $B_k$ nonsingular, it follows from
(\ref{eqn-beta1}) that $u_k$ will be equal to the right-hand side
vector up to some arbitrary non-zero scaling. Let
\begin{equation}\label{eqn-uk}
u_k= \gamma_kg_k-g_{k-1}.
\end{equation}
The scaling of $u_k$ will be reflected in $\beta_k$ by
insertion into (\ref{eqn-beta1}) as
\[
\beta_k
(\gamma_kg_k-g_{k-1})^T( g_k + \frac{g_k^T g_k}{p_{k-1}^T
  g_{k-1}}p_{k-1})= \beta_k (\gamma_k-1) g_k^T g_k =
- \frac{g_k^T g_k}{p_{k-1}^T g_{k-1}},
\]
so that
\begin{equation}\label{eqn-betak}
\beta_k=-\frac1{
(\gamma_k-1) p_{k-1}^T g_{k-1}}.
\end{equation}
Note that (\ref{eqn-betak}) is well defined as $\gamma_k\ne 1$ is
assumed. A combination of (\ref{eqn-beta1}), (\ref{eqn-uk}) and
(\ref{eqn-betak}) gives $B_k$ expressed as in (\ref{Urank1}).

It remains to show that $B_k$ is nonsingular. It follows from
(\ref{Urank1}) that
\[
B_{k} = B_{k-1} \left( I -
\frac1{
(\gamma_k-1) p_{k-1}^T g_{k-1}}
(\gamma_kB_{k-1}\inv g_k - B_{k-1}\inv g_{k-1})(\gamma_kg_k-g_{k-1})^T\right),
\]
so that
\begin{equation}\label{eqn-detBk}
\det(B_k)=\det(B_{k-1})\eta_k,
\end{equation}
with
\begin{align}\label{eqn-etak}
\eta_k 
 & = 1 -
\frac{(\gamma_kB_{k-1}\inv g_k - B_{k-1}\inv
  g_{k-1})^T(\gamma_kg_k-g_{k-1})}
{(\gamma_k-1) p_{k-1}^T g_{k-1}} \nonumber \\
& = 1 -
\frac{\gamma_k^2 g_k^T g_k  - p_{k-1}^T g_{k-1}}
{(\gamma_k-1) p_{k-1}^T g_{k-1}} 
= \frac{-\gamma_k ( \gamma_k g_k^T g_k - p_{k-1}^T g_{k-1} )}
{(\gamma_k-1) p_{k-1}^T g_{k-1}} \nonumber \\
& = \frac{-\gamma_k ( \gamma_k  - \hat\gamma_k )}
{(\gamma_k-1) \hat \gamma_k},
\end{align}
since $B_{k-1}g_k=g_k$, $B_{k-1}p_{k-1}=-g_{k-1}$ and $g_k^T
p_{k-1}=0$, with $\hat\gamma_k$ given by
Lemma~\ref{lem-gammadelta}. Hence, since $B_{k-1}$ is assumed
nonsingular, a combination of (\ref{eqn-detBk}) and (\ref{eqn-etak})
shows that nonsingularity of $B_k$ is equivalent to $\eta_k\ne 0$,
i.e., $\gamma_k\ne 0$ and $\gamma_k\ne\hat\gamma_k$, which is exactly
what is assumed.

To prove the result on positive definiteness, assume that
$B_{k-1}=B_{k-1}^T\succ0$. In this case, since $B_{k-1}$ and $B_k$
differ by a symmetric rank-1 matrix, $B_k$ can have at most one
nonpositive eigenvalue, see, e.g., \cite[Theorem
8.1.8]{golubvanloan}. Therefore, (\ref{eqn-detBk}) shows that positive
definiteness of $B_k$ is equivalent to $\eta_k>0$. Note that
$B_{k-1}\succ0$ implies $p_{k-1}^T g_{k-1}=-p_{k-1}^T B_{k-1}\inv
p_{k-1}<0$, which in turn gives $\hat\gamma_k<0$. We may now examine
(\ref{eqn-etak}) to see what values of $\gamma_k$ that give
$\eta_k>0$. The numerator of (\ref{eqn-etak}) is positive for
$\hat\gamma_k<\gamma_k<0$ and negative for $\gamma_k<\hat\gamma_k$ and
$\gamma_k>0$. The denominator of (\ref{eqn-etak}) is positive for
$\gamma_k<1$ and negative for $\gamma_k>1$.  We conclude that
$\eta_k>0$ if and only if $\hat\gamma_k<\gamma_k<0$ or $\gamma_k>1$,
which by Lemma~\ref{lem-gammadelta} is equivalent to
$0<\delta_k<\delta_k(1)$ or $\delta_k>1$.

To prove the final hereditary result, assume that $B_i p_i=-g_i$,
$i=0,\dots,k$, with $B_0=I$ and assume that $B_{i-1}$ is updated to
$B_i$ according to $(\ref{Urank1})$ for $\gamma_i$ such that
$\gamma_i\ne0$, $\gamma_i\ne\hat\gamma_i$ and $\gamma_i\ne1$.  Then,
for a given $i$, $0 < i < k$, $k$ may be replaced by $i+1$ in
(\ref{Urank1}), which gives
\begin{eqnarray}\label{eqn-Bipi}
B_{i+1}p_i & = & B_{i}p_i -
\frac1{
(\gamma_{i+1}-1) p_{i}^T g_{i}}
(\gamma_{i+1}g_{i+1} - g_{i})(\gamma_{i+1}g_{i+1}-g_{i})^Tp_i
\nonumber \\
& = & -g_i + \frac1{\gamma_{i+1}-1}(\gamma_{i+1}g_{i+1} - g_{i}) =
\frac{\gamma_{i+1}}{\gamma_{i+1}-1}(g_{i+1} - g_{i}) \nonumber \\
& = & \frac{\gamma_{i+1}\theta_i}{\gamma_{i+1}-1} H p_i,
\end{eqnarray}
where the identities $B_i p_i = - g_i$, $g_{i+1}^T p_i=0$ and
$g_{i+1}-g_i=\theta_i H p_i$ have been used. Finally, (\ref{Urank1})
gives $B_j p_i = B_{i+1} p_i$ for $j=i+2,\dots,k$, since $g_j^T p_i=0$
for $j\ge i+1$. Consequently, $B_k p_i = B_{i+1}p_i$, with $B_{i+1}
p_i$ given by (\ref{eqn-Bipi}), proving (\ref{eqn-Bkpi}). 
 \qed \end{proof}

 Note that there are two ways in which positive definiteness of a
 symmetric $B_{k-1}$ may be preserved in a symmetric rank-one
 update. The first one, $\gamma_k>1$, or equivalently
 $0<\delta_k<\delta_k(1)$, is straightforward, since it corresponds to
 $U_k\succeq 0$. The second one, $\hat\gamma_k<\gamma_k<0$, or
 equivalently $\delta_k>1$, is less straightforward. The corresponding
 $U_k$ is negative semidefinite, but still the resulting $B_k$ is
 positive definite.

 Proposition~\ref{prop-rank1} gives precise conditions for which
 rank-one matrices that give a corresponding update matrix that
 preserves positive definiteness and gives search directions parallel
 to the method of conjugate gradients. We have the freedom to choose
 $\gamma_k$ or $\delta_{k}$ appropriately. This can be compared to
 SR1, the symmetric rank-one update scheme uniquely defined by the
 secant condition
\begin{equation}\label{secantcond}
B_k s_{k-1}=y_k, \text{for} s_{k-1}=\theta_{k-1} p_{k-1}, \quad
y_k=g_k-g_{k-1}.
\end{equation}
By writing $U_k=B_k-B_{k-1}$, the secant condition gives a requirement
on $U_k$ as
\begin{equation}\label{secantcond2}
U_{k} s_{k-1}= y_k - B_{k-1}s_{k-1},
\end{equation}
which for $U_k$ symmetric and of rank one gives the SR1 update matrix
$U_k^{SR1}$ on the form
\begin{equation}\label{eqn-SR1-1}
U_k^{SR1} = \frac1{s_{k-1}^T (y_k-B_ks_{k-1})}(y_k-B_ks_{k-1})(y_k-B_ks_{k-1})^T,
\end{equation}
see, e.g., \cite[Chapter 9]{luenberger}. Since
$B_{k-1}p_{k-1}=-g_{k-1}$ holds by the definition of the quasi-Newton
method, we may use the definitions of $s_{k-1}$ and $y_k$ of
(\ref{secantcond}) to rewrite $U_k^{SR1}$ of (\ref{eqn-SR1-1}) as
\[
U_k^{SR1} =\frac1
{\theta_{k-1} p_{k-1}^T ( g_k - (1-\theta_{k-1}) g_{k-1}}
(g_k - (1-\theta_{k-1}) g_{k-1})(g_k - (1-\theta_{k-1}) g_{k-1})^T.
\]
Since exact linesearch is performed in our case, it holds that
$p_{k-1}^T g_k=0$, so that $U_k^{SR1}$ takes the form
\begin{align}\label{eqn-SR1-2}
U_k^{SR1} &=\frac{-1}
{\theta_{k-1} (1-\theta_{k-1}) p_{k-1}^Tg_{k-1}}
(g_k - (1-\theta_{k-1}) g_{k-1})(g_k - (1-\theta_{k-1}) g_{k-1})^T
\nonumber \\
&= 
\frac{-(1-\theta_{k-1})}
{\theta_{k-1} p_{k-1}^T g_{k-1}}
\left(\frac1{1-\theta_{k-1}}g_k - g_{k-1}\right)
\left(\frac1{1-\theta_{k-1}}g_k-g_{k-1}\right)^T,
\end{align}
where in the last step, a scaling of the rank-1 vector by a factor
$1/(1-\theta_{k-1})$ has been made. A comparison of (\ref{Urank1}) and
(\ref{eqn-SR1-2}) shows that the SR1 update is the particular member
of the family of symmetric rank-1 updates given by
Proposition~\ref{prop-rank1} for which
$\gamma_k=1/(1-\theta_{k-1})$. In particular, for $\theta_{k-1}=1$,
SR1 is not well defined. In addition, as there is no freedom in
choosing the rank-one matrix for SR1, there is no way to ensure $B_k
\succ 0$ even if $B_{k-1}=B_{k-1}^T\succ 0$. Note that the condition
on $B_k$ of (\ref{secantcond}) giving a condition on $U_k$ of
(\ref{secantcond2}) and a unique symmetric rank-1 $U_k$ of
(\ref{eqn-SR1-1}) is analogous to our condition on $B_k$ of
Proposition~\ref{prop-iffB} for a fixed $\delta_k$ giving a condition
on $U_k$ of Proposition~\ref{prop-iffU} and a unique rank-1 $U_k$ of
Proposition~\ref{prop-rank1}.

Example~\ref{ex-rank1} illustrates the SR1 update and another rank-1
update of Proposition~\ref{prop-rank1} which preserves positive
definiteness. The $H$ and $c$ of the example are parameterized by a
positive scalar $\phi$. We obtain $\theta_0=2/(3\phi)$, so by
selecting $\phi=2/3$, it follows that $\theta_0=1$ and the SR1 update
becomes undefined. By selecting $\phi$ slightly smaller than 2/3, for
example 0.65, we obtain $\theta_0$ slightly larger than one
($\theta_0=40/39$), so that $\gamma_1=-39$ and the corresponding
$\delta_1$ is negative ($\delta_1=-3/10$). Consequently, $B_1^{SR1}$
is indefinite and the corresponding $p_1$ is an ascent direction. For
comparison, the rank-1 update of Proposition~\ref{prop-rank1} is given
for $\delta_1=2$, which preserves positive definiteness. As can be
seen from (\ref{Urank1}), the rank-1 update of
Proposition~\ref{prop-rank1} is independent of $\phi$.

\begin{example}\label{ex-rank1}
For a positive parameter $\phi$, consider the example
\[
H= \phi \mtx{cc}{2 & \ 0 \\ 0 & \ 1}, 
\quad c= \phi \mtx{c}{-1 \\ -1}, 
\quad B_0= \mtx{cc}{1 & \ 0 \\ 0 & \ 1}, 
\quad x_0=\mtx{c}{0 \\ 0}
\]
for which
\[
\quad x_1=\mtx{r}{\frac23 \\ \frac23},
\quad x_2=\mtx{r}{\frac12 \\ 1}.
\]
Then
\begin{align*}
\phi&=\frac23 & \implies & &B_1^{SR1} &\ \mbox{undefined} &
B_1^{\delta_1=2} &=
\frac1{44}\mtx{rr}{ 43 & 5 \\ 5 & \ 19 } \\
\phi&=\frac{65}{100} & \implies & & B_1^{SR1} &=\frac1{20}
\mtx{rr}{ -16 & 42 \\ 42 & \ -29 } &
B_1^{\delta_1=2} &=
\frac1{44}\mtx{rr}{ 43 & 5 \\ 5 & \ 19 }
\end{align*}
\end{example}

Note that the numerical values and the dimension of
Example~\ref{ex-rank1} are not important. For a given quadratic
problem, there will always exist a particular positive scaling such
that the resulting $B_1^{SR1}$ is undefined. 

 \section{On the approximation of the Hessian}\label{sec-approx}

 The results of the present manuscript have been written based on the
 search directions of the method of conjugate gradients. The reason
 for doing so is that it allows a direct treatment of $B_k$, and there
 is no need to focus on the update matrix $B_k-B_{k-1}$. This is the
 choice of the authors, but other choices are of course possible.

 The results are stated for a matrix $B_k$ that approximates the
 Hessian $H$. We prefer to think of the quasi-Newton method in this
 way, but there would be little difference if one instead stated the
 results for a matrix $N_k$ that approximates $H\inv$, which is done
 for example in \cite{KON98}. The search direction $p_k$ would then be
 defined by $p_k=-N_k g_k$ rather than by $B_k p_k=-g_k$ and
 conditions would be imposed on $N_k$ rather than on
 $B_k$. Proposition~\ref{prop-iffB} could be equivalently stated using
 $N_k$ as the approximation of $H\inv$. Then, the counterparts of
 (\ref{eqn-iffB}) and (\ref{eqn-iffW}) would read
\[
N_k g_k = \delta_k A_k\inv g_k
\text{and}
A_k N_k A_k^T g_k = \delta_k g_k.
\]
When considering update matrices, with $V_k=N_k-N_{k-1}$, the
counterpart of the update formula (\ref{condU}) of
Proposition~\ref{prop-iffU} would read
\begin{align}\label{eqn-Vk}
V_k g_k &=\delta_k A_k\inv g_k - N_{k-1} g_k 
=(\delta_k-1) g_k + \delta_k \frac{g_k^T g_k}{g_{k-1}^T p_{k-1}}p_{k-1}.
\end{align}
For the rank-one case, the update matrix is unique for a given
$\delta_k$, and (\ref{eqn-Vk}) gives
\begin{equation}\label{eqn-Vkrank1}
N_k = N_{k-1} -\frac1{\gamma_k ( \gamma_k g_k^T g_k - p_{k-1}^T g_{k-1} )} (\gamma_k g_k +
  p_{k-1})(\gamma_kg_k + p_{k-1})^T,
\end{equation}
with $\gamma_k=\gamma_k(\delta_k)$ of Lemma~\ref{lem-gammadelta}. The
uniqueness of the update implies that if $N_{k-1}=B_{k-1}\inv$, then
$N_{k}=B_{k}\inv$ and (\ref{eqn-Vkrank1}) follows from (\ref{Urank1})
by the Sherman-Morrison formula.

As for the rank-one case, in light of the results of
Section~\ref{sec-rank1}, one of the referees has pointed out that the
parametrization given by $\delta_k$ can be replaced by a different
parametrization. The secant condition $\theta_{k-1} B_k p_{k-1} =
g_k-g_{k-1}$ and its counterpart on the inverse $\theta_{k-1} p_{k-1}
= N_k (g_k-g_{k-1})$ may be relaxed by a parameter $\lambda_k$ so that
$\lambda_k\theta_{k-1} B_k p_{k-1} = g_k-g_{k-1}$ and
$\lambda_k\theta_{k-1} p_{k-1} = N_k (g_k-g_{k-1})$ respectively. For
the updates, we obtain
\begin{subequations}
\begin{align}
\lambda_k \theta_{k-1} U_k p_k &= g_k
-(1-\lambda_k\theta_{k-1})g_{k-1} \text{and} \label{eqn-Ulambda} \\
V_k (g_k-g_{k-1}) &= -g_k-(1-\lambda_k\theta_{k-1})p_{k-1}, \label{eqn-Vlambda}
\end{align}
\end{subequations}
since $B_{k-1}p_{k-1}=-g_{k-1}$, $N_{k-1} g_k=g_k$ and $N_{k-1}
g_{k-1}=-p_{k-1}$. If having read the previous sections of this paper,
we would see that $p_k=\delta_k p_k^{CG}$, where we can relate
$\gamma_k$ to $\lambda_k$ by $\gamma_k = 1/(1-\lambda_k\theta_{k-1})$,
by comparing the right-hand side vector of (\ref{eqn-Ulambda}) to the
rank-one vector of Proposition~\ref{prop-rank1} or comparing the
right-hand side vector of (\ref{eqn-Vlambda}) to the rank-one vector
of (\ref{eqn-Vkrank1}). The corresponding relationship to $\delta_k$
is given by Lemma~\ref{lem-gammadelta}. An alternative to reading the
previous sections of this paper, however, would be to say that
$\lambda_k=1$ corresponds to the SR1 update, and $\lambda_k=0$
corresponds to the conjugate projection update \cite[Equation
(4.1.10)]{fletcherpractical}. They are considered in the update of the
inverse and are both known to give $p_k$ parallel to $p_k^{CG}$ if
$N_0=I$. By replacing $\gamma_k$ by $1/(1-\lambda_k\theta_{k-1})$, one
could show that a rank-1 matrix of the form (\ref{eqn-Uk}) would give
$p_k$ parallel to $p_k^{CG}$ using induction similar to what is done
in \cite[Theorem 3.4.1]{fletcherpractical} and give conditions on
preserving positive definiteness on $\lambda_k$. This would, however,
not show that there is no other family of rank-1 updates giving $p_k$
parallel to $p_k^{CG}$. We prefer to give a direct proof based on our
result of Proposition~\ref{prop-iffB}, as we from there get both
necessary and sufficient conditions.

\section{Preconditioning}\label{sec-precond}

Our results have been derived in the setting of CG, which corresponds
to $B_0=I$ in QN giving the initial search directions identical. In
this section, we give the analogous results in a preconditioned
setting. In the preconditioned method of conjugate gradients, there is
a positive definite symmetric matrix $M$, providing an estimate of
$H$. For the quasi-Newton method, this will correspond to $B_0=M$
giving the initial search directions identical.

The preconditioned method of conjugate gradient takes the following
form. If the Cholesky factor of $M$ is denoted by $L$, so that
$M=LL^T$, then the method of conjugate gradients is applied to
\begin{equation}\label{eqn-equivsystem}
L\inv H L\Tinv \xhat + L\inv c = 0,
\end{equation}
for $\xhat=L\T x$, see, e.g.,~\cite[Chapter 9.2]{saad}. Letting
``hat'' be associated with quantities of (\ref{eqn-equivsystem}), we
obtain $\phat=L\T p$ and $\ghat=L\inv g$. Since $\phat$ is associated
with a ``usual'' unpreconditioned system, we write $\phat^{CG}$, and
since $p$ is associated with a preconditioned system, we write
$p^{PCG}$, so that $\phat^{CG}=L\T p^{PCG}$. It is straightforward to
use these relations to derive the result analogous to those given in
the the previous sections also for the preconditioned system.

\begin{definition}[The preconditioned method of conjugate gradients
  (PCG)]\label{def-pcg} 
For a positive definite symmetric $n\times n$ matrix $M$, the
  \emph{preconditioned method of conjugate gradients}, PCG, is the
  linesearch method of the form given by
  Algorithm~$\ref{alg-linesearch}$ in which the search direction $p_k$
  is given by $p_k^{PCG}$, with
\begin{equation}\tag{PCG} 
p_k^{PCG}=\begin{cases}
-M\inv g_0 & \mbox{if } k=0, \\
-M\inv g_k+\frac{g_k^TM\inv g_k}{g_{k-1}^TM\inv g_{k-1}} p_{k-1}^{PCG}
& \mbox{if } k\ge 1.
\end{cases}
\end{equation}
\end{definition}
For PCG it holds that, for all $k$, $g_k^T M\inv g_i=0$, $i=0, \dots,
k-1$, so the method terminates with $g_r=0$ for some $r$, $r\le n$,
and $x_r$ solves (\ref{qp}).  In addition, it holds that
$\{p_k^{PCG}\}_{k=0}^{r-1}$ are mutually conjugate with respect to
$H$.

\begin{proposition}\label{prop-iffBprecon}
  Consider iteration $k$ of the exact linesearch method of
  Algorithm~$\ref{alg-linesearch}$, where $1\le k < r$. Assume that
  $p_i=\delta_i p_i^{PCG}$ with $\delta_i\ne 0$ for $i=0,\dots,k-1$,
  where $p_i^{PCG}$, $i=0,\dots,k-1$, are the search directions of the
  preconditioned method of conjugate gradients, as stated in
  Definition~$\ref{def-pcg}$. Let $A_k$ be
  defined as
\begin{equation*}
A_k = I-\frac{1}{g_{k-1}^T p_{k-1}} p_{k-1} g_k^T.
\end{equation*}
Then,
\begin{equation*}
A_k\inv = I+\frac{1}{g_{k-1}^T p_{k-1}} p_{k-1} g_k^T,
\end{equation*}
and it holds that $M A_k p_k^{PCG}=- g_k$. In addition, if $p_k$ is
given by $B_k p_k=-g_k$ with $B_k$ nonsingular, then, for any nonzero
scalar $\delta_k$, it holds that $p_k=\delta_k p_k^{PCG}$ if and only
if
\begin{equation*}
  B_k A_k\inv M\inv g_k = 
\frac1{\delta_k} g_k,
\end{equation*}
or equivalently if and only if
\begin{equation*}
B_k = A_k^T W_k A_k, \text{with} W_k M\inv g_k = \frac1{\delta_k} g_k,
\text{for $W_k$ nonsingular.} 
\end{equation*}
Finally, it holds that $B_k\succ 0$ if and only if $W_k\succ 0$.
\end{proposition}

In particular, $B_k = A_k^T M A_k$, corresponding to $W_k=M$ in Proposition~\ref{prop-iffBprecon}, is a
positive-definite symmetric matrix for which $B_k p_k=-g_k$ gives
$p_k=p_k^{PCG}$.

\begin{proposition}
  Consider iteration $k$ of the exact linesearch method of
  Algorithm~$\ref{alg-linesearch}$, where $1\le k < r$. Assume that
  $p_i=\delta_i p_i^{PCG}$ with $\delta_i\ne 0$ for $i=0,\dots,k-1$,
  where $p_i^{PCG}$, $i=0,\dots,k-1$, are the search directions of the
  preconditioned method of conjugate gradients using a positive
  definite symmetric preconditioning matrix $M$, as stated in
  Definition~$\ref{def-pcg}$. Let $B_{k-1}$ be a nonsingular matrix
  such that $B_{k-1} p_{k-1}=-g_{k-1}$ and $B_{k-1}M\inv g_k = g_k$. Let
  $U_{k}=B_k - B_{k-1}$ and assume that $B_k$ and $p_k$ satisfy
  $B_kp_k=-g_k$, with $B_k$ nonsingular. Then, for any nonzero scalar
  $\delta_k$, it holds that $p_k=\delta_k p_k^{PCG}$ if and only if
\begin{equation*}
U_{k} \big( M\inv g_k+\frac{g_k^TM\inv g_k}{p_{k-1}^T
  g_{k-1}}p_{k-1} \big) = (\frac1{\delta_k} -1 ) g_k + \frac{ g_k^T M\inv
  g_k}{p_{k-1}^T g_{k-1}} g_{k-1}.
\end{equation*}
\end{proposition}

\begin{lemma}\label{lem-gammadeltaprecon}
  Consider iteration $k$ of the exact linesearch method of
  Algorithm~$\ref{alg-linesearch}$, where $1\le k < r$. Assume that
  $p_i=\delta_i p_i^{PCG}$ with $\delta_i\ne 0$ for $i=0,\dots,k-1$,
  where $p_i^{PCG}$, $i=0,\dots,k-1$, are the search directions of the
  preconditioned method of conjugate gradients using a positive
  definite symmetric preconditioning matrix $M$, as stated in
  Definition~$\ref{def-pcg}$. 
Let $\hat\gamma_k={p_{k-1}^T
    g_{k-1}}/{g_k^T M\inv g_k}$. For $\delta_k\ne 0$ and
  $\gamma_k\ne\hat\gamma_k$, let the functions $\gamma_k(\delta_k)$
  and $\delta_k(\gamma_k)$ be defined by
\[
\gamma_k(\delta_k) =- \frac{p_{k-1}^T g_{k-1}}{g_k^T M\inv g_k}
\left(\frac1{\delta_k}-1\right), \quad
\delta_k(\gamma_k) =\frac1{1 - \gamma_k\frac{g_k^T M\inv
    g_k}{p_{k-1}^T g_{k-1}}}.
\]
Then, the functions $\gamma_k(\cdot)$ and $\delta_k(\cdot)$ are
inverses to each other.
\end{lemma}

\begin{proposition}\label{prop-rank1precon}
  Consider iteration $k$ of the exact linesearch method of
  Algorithm~$\ref{alg-linesearch}$, where $1\le k < r$. Assume that
  $p_i=\delta_i p_i^{PCG}$ with $\delta_i\ne 0$ for $i=0,\dots,k-1$,
  where $p_i^{PCG}$, $i=0,\dots,k-1$, are the search directions of the
  preconditioned method of conjugate gradients using a positive
  definite symmetric preconditioning matrix $M$, as stated in
  Definition~$\ref{def-pcg}$. Let
  $B_k$ and $p_k$ satisfy $B_kp_k=-g_k$, and let $B_{k-1}$ be a
  nonsingular matrix such that $B_{k-1} p_{k-1}=-g_{k-1}$ and
  $B_{k-1}M\inv g_k = g_k$. In addition, let $\gamma_k(\cdot)$,
  $\delta_k(\cdot)$ and $\hat\gamma_k$ be given by
  Lemma~$\ref{lem-gammadeltaprecon}$.

  For any scalar $\gamma_k$, except $\gamma_k=0$,
  $\gamma_k=\hat\gamma_k$ and $\gamma_k=1$, let $B_k$ be defined by
\begin{equation}\label{eqn-Urank1precon}
B_{k} = B_{k-1} -
\frac1{
(\gamma_k-1) p_{k-1}^T g_{k-1}}
(\gamma_{k}g_k - g_{k-1})(\gamma_{k}g_k-g_{k-1})^T.
\end{equation}
Then, $B_k$ is nonsingular and $p_k=\delta_k p_k^{PCG}$ for
$\delta_k=\delta_k(\gamma_k)$.

Conversely, for any scalar $\delta_k$, except $\delta_k=0$,
$\delta_k=\delta_k(1)$ and $\delta_k=1$, assume that $p_k=\delta_k
p_k^{PCG}$ and assume that $B_k-B_{k-1}$ is symmetric and of rank
one. Then, $B_k$ is a nonsingular matrix given by
$(\ref{eqn-Urank1precon})$ for $\gamma_k=\gamma_k(\delta_k)$.

If, in addition, $B_{k-1}=B_{k-1}^T\succ0$, then $B_k$ defined by
$(\ref{eqn-Urank1precon})$ satisfies $B_k\succ 0$ if and only if
$\gamma_k>1$ or $\hat\gamma_k<\gamma_k<0$, or equivalently if and only
if $\gamma_k=\gamma_k(\delta_k)$ for $0<\delta_k<\delta_k(1)$ or
$\delta_k>1$.

Finally, if $B_i p_i=-g_i$, $i=0,\dots,k$, with $B_0=M$ and if, for
$i=1,\dots,k$, $B_{i-1}$ is updated to $B_i$ according to
$(\ref{eqn-Urank1precon})$ for $\gamma_i$ such that $\gamma_i\ne0$,
$\gamma_i\ne\hat\gamma_i$ and $\gamma_i\ne1$, then
\[
B_k p_i = \frac{\gamma_{i+1}\theta_i}{\gamma_{i+1}-1}H p_i, \quad
i=0,\dots,k-1.
\]
\end{proposition}

\section{Conclusion}\label{Conclusion}

In this paper we have derived necessary and sufficient conditions on
the matrix $B_k$ in a QN-method such that $p_k$, obtained by solving
$B_k p_k=-g_k$, satisfies $p_k=\delta_k p_k^{PCG}$ for some $\delta_k
\neq 0$, where $p_k^{PCG}$ is the search direction of the
preconditioned method of conjugate gradients. These conditions are stated in
Proposition~\ref{prop-iffBprecon}. The results have been derived for
the case of CG and then extended to PCG for a symmetric positive
definite preconditioning matrix $M$.

Further, we have characterized the symmetric rank-one update matrices
for QN that give parallel search directions to those of PCG. In
Proposition~\ref{prop-rank1precon}, we show that the rank-one matrix
must be a linear combination of $g_k$ and $g_{k-1}$, and also that
almost any linear combination will do. In addition, we characterize
the family of symmetric rank-one updates that preserve symmetry and
positive definiteness of $B_{k-1}$.

Our focus is on the mathematical properties of PCG and QN in exact
arithmetic. We want to stress that considering the numerical
properties in finite precision is of utmost importance, but such an
analysis is beyond the scope of this paper.  See, e.g., \cite{hager}
for an illustration of a case where PCG and QN generate identical
iterates in exact arithmetic but the difference between numerically
computed iterates for the two methods is large.

The results of the paper are meant to be useful as such, for
understanding the behavior of exact linesearch quasi-Newton methods
for minimizing a quadratic function. In addition, we hope that they
can lead to further research on methods for unconstrained
minimization. In particular, understanding the behavior of
quasi-Newton methods on near-quadratic functions would be a subject of
future research.

\section*{Acknowledgements}

We thank the editor and the anonymous referees for their constructive
comments which significantly improved the presentation.


\end{document}